\documentclass[11pt]{article}

\usepackage[centertags]{amsmath}
\usepackage{amsfonts}
\usepackage{amssymb}
\usepackage{amsthm}
\usepackage{newlfont}

\newtheorem{thm}{Theorem}[section]
\newtheorem*{thm*}{Theorem}
\newtheorem{cor}[thm]{Corollary}

\newtheorem{lem}[thm]{Lemma}

\newtheorem{prop}[thm]{Proposition}
\newtheorem*{prop*}{Proposition}

\newtheorem*{conj*}{Conjecture}

\newtheorem*{dfn*}{Definition}
\theoremstyle{definition}
\newtheorem{rem}[thm]{\textbf{Remark}}
\newtheorem*{rmk*}{Remark}
\newtheorem*{fact*}{Fact}

\theoremstyle{proof}

\newcommand{\norm}[1]{\left\Vert#1\right\Vert}

\newcommand{\abs}[1]{\left\vert#1\right\vert}
\newcommand{\set}[1]{\left\{#1\right\}}
\newcommand{\brac}[1]{\left(#1\right)}

\newcommand{\Real}{\mathbb{R}}

\newcommand{\eps}{\varepsilon}

\newcommand{\I}{\mathcal{I}}

\newcommand{\B}{\mathcal{B}}

\def \F {\mathcal{F}}
\def \J {\mathcal{N}}

\newlength{\defbaselineskip}
\newcommand{\setlinespacing}[1]%
           {\setlength{\baselineskip}{#1 \defbaselineskip}}

\numberwithin{equation}{section}

\begin{document}

\title{A converse to Maz'ya's inequality for capacities under curvature lower bound}
\author{Emanuel Milman\textsuperscript{1} \\ \\
\emph{Dedicated to V.G. Maz'ya with great respect and admiration}}

\footnotetext[1]{School of Mathematics,
Institute for Advanced Study, Einstein Drive, Simonyi Hall, Princeton, NJ 08540, USA.
Email: emilman@math.ias.edu.\\
Supported by NSF under agreement \#DMS-0635607. }

\maketitle

\begin{abstract}
We survey some classical inequalities due to Maz'ya relating isocapacitary inequalities with
their functional and isoperimetric counterparts in a measure-metric space setting, and extend Maz'ya's lower bound for the $q$-capacity ($q>1$) in terms of the $1$-capacity (or isoperimetric) profile. We then proceed to describe results by Buser, Bakry, Ledoux and most recently by the author, which show that under suitable convexity assumptions on the measure-metric space, Maz'ya's inequality for capacities may be reversed, up to dimension independent numerical constants: a
matching lower bound on $1$-capacity may be derived in terms of the $q$-capacity profile. We extend these results to handle arbitrary $q > 1$ and weak semi-convexity assumptions, by obtaining some new delicate semi-group estimates.  
\end{abstract}


\section{Introduction}

The notion of capacity, first systematically introduced by Frechet, has played a fundamental role in the 
theory developped by V. G. Maz'ya in the 1960's for the study of functional inequalities and 
embedding theorems, and has continued to play an important role in the development of the theory ever since
(see \cite{MazyaBook} for an extended overview). 
Before recalling the definition, let us first describe our setup.

We will denote by $(\Omega,d)$ a separable metric space, and by
$\mu$ a Borel \emph{probability} measure on $(\Omega,d)$ which
is not a unit mass at a point. Let $\F = \F(\Omega,d)$ denote the space of functions which are Lipschitz on every ball in $(\Omega,d)$ - we will call such functions ``Lipschitz-on-balls''. Given $f \in \F$, we will denote by $\abs{\nabla f}$ the following Borel function:
\[
 \abs{\nabla f}(x) := \limsup_{d(y,x) \rightarrow 0+} \frac{|f(y) - f(x)|}{d(x,y)} ~.
\]
(and we define it as 0 if $x$ is an isolated point - see \cite[pp. 184,189]{BobkovHoudre}
for more details). Although it is not essential for the
ensuing discussion, it will be more convenient to think of $\Omega$ as a complete smooth oriented $n$-dimensional
Riemannian manifold $(M,g)$ and of $d$ as the induced geodesic distance, in which case $\abs{\nabla f}$ coincides with the 
usual Riemannian length of the gradient. 

\begin{dfn*}
Given two Borel sets $A \subset B \subset (\Omega,d)$ and $1 \leq q<\infty$, the $q$-capacity of $A$ relative to $B$ is defined as:
\[
 Cap_q(A,B) := \inf\set{ \norm{\abs{\nabla \Phi}}_{L_q(\mu)} ; \Phi|_A \equiv 1 \;,\;
\Phi|_{\Omega \setminus B} \equiv 0 },
\]
where the infimum is on all $\Phi : \Omega \rightarrow [0,1]$ which are Lipschitz-on-balls.
\end{dfn*}

We remark that it is possible to give an even more general definition than the one above (see Maz'ya \cite{MazyaBook}).
Note that in the case of a compact manifold $(M,g)$ and the Riemmanian (normalized) volume $\mu$, $Cap_2(A,B)^2$ coincides (up to constants) with the usual Newtonian capacity of a compact set $A$ relative to the outer open set $B$.  Following Barthe--Cattiaux--Roberto \cite{BCRHard}, we will only be interested in this work in the \emph{$q$-capacity profile}: 

\begin{dfn*}
Given a metric probability space $(\Omega,d,\mu)$, $1 \leq q<\infty$, its $q$-capacity profile is defined for any $0 < a \leq b < 1$ as:
\begin{eqnarray*}
 Cap_q(a,b) & := & \inf\set{ Cap_q(A,B) \; ; \; A \subset B \; , \; \mu(A) \geq a \;,\; \mu(B) \leq b } \\
 & = & \inf\set{ \norm{\abs{\nabla \Phi}}_{L_q(\mu)} \; ; \; \mu\set{\Phi=1} \geq a \;,\; \mu\set{\Phi=0} \geq 1-b },
\end{eqnarray*}
where the latter infimum is on all $\Phi : \Omega \rightarrow [0,1]$ which are Lipschitz-on-balls.
\end{dfn*}

The intimate relation between $1$-capacity and the isoperimetric properties of a space was noticed by Fleming \cite{FlemingBVFunctions} in the Euclidean setting, using the co-area formula of Federer \cite{FedererCurvatureMeasures} (see also Federer--Fleming \cite{FedererFleming}), and generalized by Maz'ya \cite{MazyaSobolevImbedding}. An analogous relation between isoperimetric inequalities and functional inequalities involving the term $\norm{\abs{\nabla f}}_{L_1(\mu)}$ was discovered by Maz'ya \cite{MazyaSobolevImbedding} and independently by Federer--Fleming \cite{FedererFleming}, leading in particular to the determination of the optimal constant in the Gagliardo inequality in Euclidean space $(\Real^n,\abs{\cdot})$. Maz'ya continued to study these relations when $1$ above is replaced by a general $q>1$ \cite{MazyaSobolevImbedding,MazyaCapacities,MazyaCheegersInq1,MazyaBook}: he showed how to pass from any lower bound on $1$-capacity to an optimal lower bound on $q$-capacity, and demonstrated the equivalence between $q$-capacitary inequalities and functional inequalities involving the term $\norm{\abs{\nabla f}}_{L_q(\mu)}$. Combining all these ingredients, Maz'ya discovered a way to pass from isoperimetric information to optimal (up to constants) functional inequalities. Especially useful is the case $q=2$, since this corresponds to spectral information on the Laplacian and the Schr\"{o}dinger operators, leading to many classical characterizations \cite{MazyaCheegersInq1,MazyaCheegersInq2,MazyaBook}.

We will define all of the above notions in Section \ref{sec:defs}, and sketch the proofs of their various relations in our  metric probability space setting in Section \ref{sec:caps}. Moreover, we will extend the the transition from $1$-capacity to $q$-capacity ($q>1$) to handle arbitrary transition between $p$ and $q$, when $p<q$, following \cite{EMilmanRoleOfConvexityInFunctionalInqs}. 

\medskip

It is easy to check that in general, one cannot deduce back information on $p$-capacity from $q$-capacity, when $p<q$. In other words, the above transition is only one-directional, and cannot in general be reversed (see Subsection \ref{subsec:known-connections}). We will therefore need to add some additional assumptions in order to have any chance of obtaining a reverse implication. As we will see below, some type of convexity assumptions are a natural candidate. We start with two important examples when $(M,g) =
(\Real^n,\abs{\cdot})$ and $\abs{\cdot}$ is some fixed Euclidean norm:
\begin{itemize}
\item
$\Omega$ is an \emph{arbitrary} bounded convex domain in $\Real^n$
($n \geq 2$), and $\mu$ is the uniform probability measure on $\Omega$.
\item
$\Omega = \Real^n$ ($n \geq 1$) and $\mu$ is an \emph{arbitrary}
absolutely continuous log-concave probability measure, meaning that
$d\mu = \exp(-\psi) dx$ where $\psi: \Real^n \rightarrow \Real \cup
\set{+\infty}$ is convex (we refer to the paper
\cite{Borell-logconcave} of C. Borell for more information).
\end{itemize}

In both cases, we will say that ``our convexity assumptions are
fulfilled''. More generally, we will use the following definition from
\cite{EMilmanGeometricApproachPartI}:

\begin{dfn*}
We will say that our \emph{smooth $\kappa$-semi-convexity assumptions} are
fulfilled if:
\begin{itemize}
\item
$(M,g)$ denotes an $n$-dimensional ($n\geq 2$) oriented smooth complete connected
Riemannian manifold or $(M,g)=(\Real,\abs{\cdot})$, and $\Omega =
M$.
\item
$d$ denotes the induced geodesic distance on $(M,g)$.
\item
$d\mu = \exp(-\psi)  dvol_M$,  $\psi \in C^2(M)$, and as tensor
fields on $M$:
\begin{equation} \label{eq:Intro-BE}
Ric_g + Hess_g \psi \geq -\kappa g ~.
\end{equation}
\end{itemize}
We will say that our \emph{$\kappa$-semi-convexity assumptions} are fulfilled if
$\mu$ can be approximated in total-variation by measures
$\set{\mu_m}$ so that each $(\Omega,d,\mu_m)$ satisfies our smooth
$\kappa$-semi-convexity assumptions. \\
When $\kappa=0$, we will say that our (smooth) convexity assumptions are satisfied. 
\end{dfn*}

The condition (\ref{eq:Intro-BE}) is the well-known
Curvature-Dimension condition $CD(-\kappa,\infty)$, introduced by Bakry
and \'Emery in their celebrated paper \cite{BakryEmery} (in the more
abstract framework of diffusion generators). Here $Ric_g$ denotes
the Ricci curvature tensor and $Hess_g$ denotes the second covariant
derivative.   

\medskip

Our main result from \cite{EMilman-RoleOfConvexity}, as extended in \cite{EMilmanRoleOfConvexityInFunctionalInqs}, is that under our convexity assumptions ($\kappa = 0$ case), the above transition can in fact be reversed, up to \emph{dimension independent} constants. This can be formulated in terms of passing from $q$-capacity to $p$-capacity ($1\leq p<q$), or equivalently, as passing from functional inequalities involving the term $\norm{\abs{\nabla f}}_{L_q(\mu)}$ to isoperimetric inequalities. 

In this work, we extend our previous results to handle the more general $\kappa$-semi-convexity assumptions. For the case $q=2$, this was previously done by Buser \cite{BuserReverseCheeger} in the case of a uniform density on a manifold with Ricci curvature bounded from below, and extended to the general smooth $\kappa$-semi-convexity assumptions by Bakry--Ledoux \cite{BakryLedoux} and Ledoux \cite{LedouxSpectralGapAndGeometry} using a diffusion semi-group approach. To handle the general $q>1$ case, we follow the semi-group argument, as in our previous work \cite{EMilmanRoleOfConvexityInFunctionalInqs}. Surprisingly, the case $\kappa > 0$ requires proving new delicate semi-group estimates, which may be of independent interest, and which were not needed for the previous arguments. We formulate and prove this converse to Maz'ya's inequality for capacities in Section \ref{sec:new}.

\medskip

\textbf{Acknowledgements.} I would like to thank Sasha Sodin for introducing me to Maz'ya's work on capacities. 

\section{Definitions and Preliminaries} \label{sec:defs}

\subsection{Isoperimetric Inequalities}

In Euclidean space, an isoperimetric inequality relates between the (appropriate notion of) surface area of a Borel set and its volume. To define an appropriate generalization of surface area in our setting, we will use Minkowski's (exterior) boundary measure of a Borel set $A \subset (\Omega,d)$, denoted here by $\mu^+(A)$, which is defined as:
\[
 \mu^+(A) := \liminf_{\eps \to 0} \frac{\mu(A^d_{\eps}) - \mu(A)}{\eps}~,
\]
where $A^d_{\eps} := \set{x \in \Omega ; \exists y \in A \;\; d(x,y)
< \eps}$ denotes the $\eps$-neighborhood of $A$ with respect to the
metric $d$. An isoperimetric inequality measures the
relation between $\mu^+(A)$ and $\mu(A)$ by means of the
isoperimetric profile $I = I_{(\Omega,d,\mu)}$, defined as the
pointwise maximal function $I : [0,1] \rightarrow \Real_+$, so that:
\begin{equation} \label{eq:I-definition}
\mu^+(A) \geq I(\mu(A)) ~,
\end{equation}
for all Borel sets $A \subset \Omega$. Since $A$ and $\Omega \setminus A$ will typically have the same boundary measure, it will be convenient to also define $\tilde{\I} :[0,1/2] \rightarrow \Real_+$ as $\tilde{\I}(v) := \min(\I(v),\I(1-v))$.

Let us keep some important examples of isoperimetric inequalities in mind. 
We will say that our space satisfies \emph{linear isoperimetric inequality}, if there exists a
constant $D>0$ so that $\tilde{I}_{(\Omega,d,\mu)}(t) \geq D t$ for
all $t \in [0,1/2]$ ; we denote the best constant $D$ by $D_{Lin} =
D_{Lin}(\Omega,d,\mu)$. Another useful example pertains to the
standard Gaussian measure $\gamma$ on $(\Real,\abs{\cdot})$, where
$\abs{\cdot}$ is the Euclidean metric. We will say that our space
satisfies a \emph{Gaussian isoperimetric inequality}, if there
exists a constant $D>0$ so that $I_{(\Omega,d,\mu)}(t) \geq D
I_{(\Real,\abs{\cdot},\gamma)}(t)$ for all $t \in [0,1]$ ; we denote
the best constant $D$ by $D_{Gau}= D_{Gau}(\Omega,d,\mu)$. It is known that
$\tilde{I}_{(\Real,\abs{\cdot},\gamma)}(t) \simeq t \log^{1/2}(1/t)$
uniformly on $t \in [0,1/2]$, where we use the notation $A \simeq B$
to signify that there exist universal constants $C_1,C_2>0$ so that
$C_1 B \leq A \leq C_2 B$. Unless otherwise stated, all of the
constants throughout this work are universal, independent of any
other parameter, and in particular the dimension $n$ in the case of
an underlying manifold. The Gaussian isoperimetric inequality can
therefore be equivalently stated as asserting that there exists a
constant $D>0$ so that $\tilde{I}_{(\Omega,d,\mu)}(t) \geq D t
\log^{1/2}(1/t)$ for all $t \in [0,1/2]$.

\subsection{Functional Inequalities}

Let $f \in \F$. We will consider functional inequalities which compare between
$\norm{f}_{N_1(\mu)}$ and $\norm{\abs{\nabla f}}_{N_2(\mu)}$, where
$N_1, N_2$ are some norms associated with the measure $\mu$, like
the $L_p(\mu)$ norms, or some other more general \emph{Orlicz
quasi-norms} associated to the class $\J$ of increasing continuous
functions mapping $\Real_+$ onto $\Real_+$.

A function $N: \Real_+ \rightarrow \Real_+$ will
be called a Young function if $N(0)=0$ and $N$ is convex increasing. Given a Young function $N$,
the Orlicz norm $N(\mu)$ associated to $N$ is defined as:
\[
 \norm{f}_{N(\mu)} := \inf \set{ v>0 ; \int_\Omega N(|f|/v) d\mu \leq 1}.
\]
For a general increasing continuous function $N : \Real_+
\rightarrow \Real_+$ with $N(0)=0$ and $\lim_{t \rightarrow \infty}
N(t) = \infty$ (we will denote this class by $\J$), the above definition still makes sense, although
$N(\mu)$ will no longer necessarily be a norm. We will say in this
case that it is a quasi-norm.

There is clearly no point to test constant functions in our functional inequalities, 
so it will be natural to require that either the expectation $E_\mu f$ or median $M_\mu f$
of $f$ are 0. Here $E_\mu f = \int f d\mu$ and $M_\mu f$ is a value
so that $\mu(f \geq M_\mu f) \geq 1/2$ and $\mu(f \leq M_\mu f) \geq
1/2$.

\begin{dfn*}
We will say that the space $(\Omega,d,\mu)$ satisfies an $(N,q)$
Orlicz-Sobolev inequality ($N \in \J, q \geq 1$) if:
\begin{equation} \label{eq:OS-inq-def}
\exists D>0 \; \text{ s.t. } \; \forall f \in \F \;\;\;\;  D \norm{f
- M_\mu f}_{N(\mu)} \leq \norm{\abs{\nabla f}}_{L_q(\mu)} ~.
\end{equation}
\end{dfn*}
\noindent A similar (yet different) definition was given by Roberto
and Zegarlinski \cite{RobertoZegarlinski} in the case $q=2$
following the work of Maz'ya \cite[p. 112]{MazyaBook}. Our preference to use the median $M_\mu$ in our definition (in place of the more standard expectation $E_\mu$) is immaterial whenever $N$
is a convex function, due to the following elementary lemma from \cite{EMilman-RoleOfConvexity}:
\begin{lem} \label{lem:E-M}
Let $N(\mu)$ denote an Orlicz norm associated to the Young function $N$. Then:
\[
 \frac{1}{2} \norm{f - E_\mu f}_{N(\mu)} \leq \norm{f - M_\mu f}_{N(\mu)} \leq 3 \norm{f - E_\mu
f}_{N(\mu)}.
\]
\end{lem}

\medskip

When $N(t)=t^p$, then $N(\mu)$ is just the usual $L_p(\mu)$
norm. If in addition $M_\mu$ in (\ref{eq:OS-inq-def}) is replaced by
$E_\mu$, the case $p=q=2$ is then just the classical \emph{Poincar\'e inequality}, and we denote the best constant in this inequality by $D_{Poin}$. Similarly, the case $q=1, p=\frac{n}{n-1}$ corresponds to the \emph{Gagliardo inequality}, 
and the case $1<q<n$, $p= \frac{qn}{n-q}$ to the \emph{Sobolev inequalities}. A limitting case when $q=2$ and $n$ tends to infinity is the so-called \emph{log-Sobolev inequality}.
More generally, we say that our space satisfies a \emph{$q$-log-Sobolev
inequality} ($q \in [1,2]$), if there exists a constant $D>0$ so
that:
\begin{equation} \label{eq:q-log-Sob}
\forall f \in \F \;\;\;\;  D \brac{\int |f|^q \log |f|^q d\mu - \int
|f|^q d\mu \log(\int |f|^q d\mu)}^{1/q} \leq \norm{\abs{\nabla
f}}_{L_q(\mu)} ~.
\end{equation}
The best possible constant $D$ above is denoted by $D_{LS_q} =
D_{LS_q}(\Omega,d,\mu)$. Although these inequalities do not
precisely fit into our announced framework, it follows from the work
of Bobkov and Zegarlinski \cite[Proposition 3.1]{BobkovZegarlinski} (generalizing the
case $q=2$ due to Bobkov and G\"{o}tze \cite[Proposition 4.1]{BobkovGotzeLogSobolev}) that they are in
fact equivalent to the following Orlicz-Sobolev inequalities:
\begin{equation} \label{eq:N_q-log-Sob}
\forall f \in \F \;\;\; D_{\varphi_q} \norm{ f - E_\mu f }_{\varphi_q(\mu)} \leq \norm{\abs{\nabla
f}}_{L_q(\mu)},
\end{equation}
where $\varphi_q(t) = t^q \log(1+t^q)$, and $D_{LS_q} \simeq D_{\varphi_q}$ uniformly on $q \in [1,2]$.

Various other functional
inequalities admit an equivalent (up to universal constants)
formulation using an appropriate Orlicz norm $N(\mu)$ on the left
hand side of (\ref{eq:OS-inq-def}). We refer the reader to the
recent paper of Barthe and Kolesnikov \cite{BartheKolesnikov} and
the references therein for an account of several other types of
functional inequalities.

\subsection{Known Connections} \label{subsec:known-connections}

It is well known that various isoperimetric inequalities imply their
functional ``counterparts''. It was shown by Maz'ya
\cite{MazyaCheegersInq1,MazyaCheegersInq2} and independently by
Cheeger \cite{CheegerInq}, that a linear isoperimetric inequality implies Poincar\'e's inequality:
$D_{Poin} \geq D_{Lin}/2$ (Cheeger's inequality). It was first
observed by M. Ledoux \cite{LedouxBusersTheorem} that a Gaussian
isoperimetric inequality implies a $2$-log-Sobolev inequality:
$D_{LS_2} \geq c D_{Gau}$, for some universal constant $c>0$. This
has been later refined by Beckner (see
\cite{LedouxLectureNotesOnDiffusion}) using an equivalent functional
form of the Gaussian isoperimetric inequality due to S. Bobkov
\cite{BobkovFunctionalFormOfGaussianIsopInq,BobkovGaussianIsopInqViaCube}
(see also  \cite{BartheMaureyIsoperimetricInqs}): $D_{LS_2} \geq
D_{Gau} / \sqrt{2}$. The constants $2$ and $\sqrt{2}$ above are
known to be optimal. 

Another example is obtained by considering the isoperimetric inequality $\tilde{I}(t) \geq D_{Exp_q} t \log^{1/q} 1/t$, for some $q \in [1,2]$. This inequality is satisfied with some $D_{Exp_q} > 0$ by the probability measure $\mu_p$ with density $\exp(-|x|^p)/Z_p$ on $(\Real,\abs{\cdot})$, for $p = q^* = q/(q-1)$ (where $Z_p$ is a normalization factor). Bobkov and Zegarlinski \cite{BobkovZegarlinski} have shown that $D_{LS_q} \geq c D_{Exp_q}$, for some universal constant $c>0$, independent of $q$, in analogy to the inequality $D_{LS_2} \geq c D_{Gau}$ mentioned above. Another proof of this using capacities was given in our joint work with Sasha Sodin \cite{EMilmanSodinIsoperimetryForULC}.   

We will see in Section \ref{sec:caps} how Maz'ya's general framework may be used to obtain all of these implications.

\medskip

In general, however, it is known that these implications \emph{cannot} be
reversed. For instance, using $([-1,1],\abs{\cdot},\mu_\alpha)$
where $d\mu_\alpha = \frac{1+\alpha}{2} |x|^\alpha dx$ on $[-1,1]$,
clearly $\mu_\alpha^+([0,1]) = 0$ so $D_{Lin} = D_{Gau} = 0$,
whereas one can show that $D_{Poin},D_{LS_2} > 0$ for $\alpha \in
(0,1)$ using criteria for the Poincar\'e and $2$-log-Sobolev
inequalities on $\Real$ due to 
Artola, Talenti and Tomaselli (cf. Muckenhoupt \cite{MuckenhouptHardyInq}) and Bobkov
and G\"{o}tze \cite{BobkovGotzeLogSobolev}, respectively. These examples suggest that we must
rule out the existence of narrow ``necks'' in our measure or space, for the converse implications to stand a chance
of being valid. Adding some convexity-type assumptions is therefore a natural path to take. 

Indeed, under our $\kappa$-semi-convexity assumptions, we will see that a reverse implication \emph{can} in fact be obtained. 
This extends some previously known results by several authors.
Buser showed in \cite{BuserReverseCheeger} that $D_{Lin} \geq c \min(D_{Poin} ,D_{Poin}^2 / \sqrt{\kappa})$ with $c>0$ a universal constant, for the case of a compact Riemannian manifold $(M,g)$ with uniform density whose Ricci curvature is bounded below by $-\kappa g$. This was subsequently extended by Ledoux \cite{LedouxSpectralGapAndGeometry} to our more general smooth $\kappa$-semi-convexity assumptions, following the semi-group approach he developped in \cite{LedouxBusersTheorem} and refined by Bakry--Ledoux \cite{BakryLedoux}. Similarly, the reverse inequality 
$D_{Gau} \geq c \min(D_{LS_2} , D_{LS_2}^2 / \sqrt{\kappa} )$ with $c>0$ a universal constant, was obtained by Bakry--Ledoux (see also \cite{LedouxSpectralGapAndGeometry}) under our smooth $\kappa$-semi-convexity assumptions. We will extend these results to handle general Orlicz-Sobolev inequalities in Section \ref{sec:new}, and in particular show that $D_{Exp_q} \geq c \min(D_{LS_q},D_{LS_q}^2 / \sqrt{\kappa} )$ for some universal constant $c>0$, uniformly on $q \in [1,2]$. When $q > 2$, we will see that these formulae take on a different form.  


\section{Capacities} \label{sec:caps}

In this section, we formulate the various known connections mentioned in the Introduction between capacities and isoperimetric and functional inequalities, and provide for completeness most of the proofs, following \cite{EMilmanRoleOfConvexityInFunctionalInqs}. 

\begin{rem} \label{rem:cap-app0}
A remark which will be useful for dealing with general metric
probability spaces, is that in the definition of capacity, by approximating $\Phi$ appropriately, we may
always assume that $\int_{\set{\Phi=t}} \abs{\nabla \Phi}^q d\mu =
0$, for any $t \in (0,1)$, even though we may have
$\mu\set{\Phi=t}>0$. See \cite[Remark 3.3]{EMilmanRoleOfConvexityInFunctionalInqs} for more information.
\end{rem}


\subsection{$1$-capacity and isoperimetric profiles}

Our starting point is the following well-known co-area formula, which in our setting becomes an inequality (see Bobkov--Houdr\'e \cite{BobkovHoudre,BobkovHoudreMemoirs}):

\begin{lem}[Bobkov--Houdr\'e]
For any $f \in \F$, we have:
\[
\int_\Omega \abs{\nabla f} d\mu \geq \int_{-\infty}^{\infty}
\mu^+\set{f>t} dt ~. 
\]
\end{lem}

The following proposition (see \cite{MazyaCapacities}, \cite{FedererFleming}, \cite{BobkovHoudre}) encapsulates the connection
between the $1$-capacity and isoperimetric profiles.

\begin{prop}[Federer--Fleming, Maz'ya, Bobkov--Houdr\'e] \label{prop:Sodin-prop}
For all $0<a < b < 1$:
\begin{equation} \label{eq:Sodin-prop}
 \inf_{a \leq t \leq b} I(t) \leq Cap_1(a,b) \leq \inf_{a \leq t < b}
I(t) ~.
\end{equation}
\end{prop}

For completeness, we provide a proof following Sodin \cite[Proposition A]{SodinLpIsoperimetry}.

\begin{proof}
Given a function $\Phi:\Omega \rightarrow [0,1]$ which is Lipschitz-on-balls
with $\mu\set{\Phi=1} \geq a$ and  $\mu\set{\Phi=0} \geq 1-b$, the
co-area inequality implies:
\[
\int \abs{\nabla \Phi} d\mu \geq \int_{-\infty}^{\infty}
\mu^+\set{\Phi>t} dt = \int_{0}^{1} \mu^+\set{\Phi>t} dt \geq
\inf_{a \leq t \leq b} I(t) ~.
\]
Taking infimum on all such functions $\Phi$, the first inequality in
(\ref{eq:Sodin-prop}) follows. To obtain the second inequality, let $A$
denote a Borel set with $a \leq \mu(A) < b$. We may exclude the case
that $\mu^+(A) = \infty$, since it does not contribute to the
definition of the isoperimetric profile $I$. Now denote for $r,s>0$:
\[
\Phi_{r,s}(x) := \brac{1 - s^{-1} d(x,A_r)} \vee 0 ~.
\]
Clearly $\mu\set{\Phi_{r,s}=1} \geq \mu(A) \geq a$, and since
$\mu^+(A) < \infty$, for $r+s$ small enough we have
$\mu\set{\Phi_{r,s}=0} \geq 1-b$. Hence:
\[
\frac{\mu(A_{s+2r}) - \mu(A)}{s} \geq \frac{\mu\set{r \leq d(x,A)
\leq s+r}}{s} \geq \int \abs{\nabla \Phi_{r,s}} d\mu \geq Cap_1(a,b)
~.
\]
Taking the limit inferior as $r,s \rightarrow 0$ so that $r/s
\rightarrow 0$, and taking infimum on all sets $A$ as above, the
second inequality in (\ref{eq:Sodin-prop}) follows.
\end{proof}

Since obviously $Cap_1(a,b) = Cap_1(1-b,1-a)$, we have the following useful corollary: 

\begin{cor} \label{cor:cap1}
For any non-decreasing continuous function $J:[0,1/2]\rightarrow \Real_+$:
\[
 \tilde{I}(t) \geq J(t) \;\;\; \forall t \in [0,1/2] \;\; \iff \;\; Cap_1(t,1/2)
\geq J(t) \;\;\; \forall t \in [0,1/2] ~.
\]
\end{cor}

\begin{dfn*}
A $q$-capacitary inequality is an inequality of the form:
\[
Cap_q(t,1/2) \geq J(t) \;\;\; \forall t \in [0,1/2] ~,
\]
where $J : [0,1/2] \rightarrow \Real_+$ is a non-decreasing continuous function.  
\end{dfn*}


\subsection{$q$-capacitary and weak Orlicz-Sobolev inequalities}

\begin{dfn*}
Given $N \in \J$, denote by $N^\wedge : \Real_+ \rightarrow \Real_+$ the ``adjoint'' function:
\[
 N^\wedge(t) := \frac{1}{N^{-1}(1/t)}.
\]
\end{dfn*}

\begin{rem} \label{rem:J-wedge}
Note that the operation $N \rightarrow N^\wedge$ is an involution on $\J$, and that
$N(\cdot^\alpha)^\wedge = (N^\wedge)^{1/\alpha}$ for $\alpha > 0$. It is also immediate to check that
$N(t^\alpha)/t$ is non-decreasing iff $N^\wedge(t)^{1/\alpha}/t$ is non-increasing ($\alpha>0$).
\end{rem}

We denote by $L_{s,\infty}(\mu)$ the weak $L_s$ quasi-norm, defined as:
\[
 \norm{f}_{L_{s,\infty}(\mu)} := \sup_{t>0} \mu(|f| \geq t)^{1/s} t.
\]
We now extend the definition of the weak $L_{s}$ quasi-norm to Orlicz quasi-norms $N(\mu)$, using the adjoint
function $N^\wedge$:
\begin{dfn*}
Given $N\in \J$, define the weak $N(\mu)$ quasi-norm as:
\[
 \norm{f}_{N(\mu),\infty} := \sup_{t>0} N^\wedge(\mu\set{\abs{f} \geq t}) t.
\]
\end{dfn*}
\noindent This definition is consistent with the one for $L_{s,\infty}$, and satisfies:
\begin{equation} \label{eq:weak-norm}
 \norm{f}_{N(\mu),\infty} \leq \norm{f}_{N(\mu)} ~,
\end{equation}
as easily checked using the Markov-Chebyshev inequality. Also note that this is indeed a quasi-norm by a simple
union-bound:
\[
\norm{f+g}_{N(\mu),\infty} \leq 2 \brac{\norm{f}_{N(\mu),\infty} + \norm{g}_{N(\mu),\infty}} ~.
\]

\begin{rem}
The motivation for the definition of $N^\wedge$ stems from the immediate observation that for any
Borel set $A$:
\[
 \norm{\chi_A}_{N(\mu)} =  \norm{\chi_A}_{N(\mu),\infty} = N^\wedge(\mu(A)) ~.
\]
For this reason, the expression $1 / N^{-1}(1/t)$ already appears in the works of Maz'ya \cite[p. 112]{MazyaBook}
and Roberto--Zegarlinski \cite{RobertoZegarlinski}.
\end{rem}

\begin{dfn*}
An inequality of the  form:
\begin{equation} \label{eq:Lq-implies-capq}
\forall f \in \F \;\;\; D \norm{f - M_\mu f}_{N(\mu),\infty} \leq
\norm{\abs{\nabla f}}_{L_q(\mu)}
\end{equation}
is called a weak-type Orlicz-Sobolev inequality.
\end{dfn*}

\begin{lem} \label{lem:Lq-implies-capq}
The weak-type Orlicz-Sobolev inequality (\ref{eq:Lq-implies-capq}) implies:
\[
 Cap_q(t,1/2) \geq D N^\wedge(t) \;\;\; \forall t \in [0,1/2]
\]
\end{lem}

\begin{proof}
Apply (\ref{eq:Lq-implies-capq}) to $f = \Phi$, where $\Phi : \Omega \rightarrow
[0,1]$ is any Lipschitz-on-balls function so that $\mu\set{\Phi=1} \geq t$ and
$\mu\set{\Phi=0} \geq 1/2$. Since $M_\mu \Phi = 0$, it
follows that:
\[
 \norm{\abs{\nabla \Phi}}_{L_q(\mu)} \geq D \norm{\Phi}_{N(\mu),\infty} \geq D
N^\wedge(\mu(\set{\Phi=1})) \geq D N^\wedge(t),
\]
Taking the infimum over all $\Phi$ as above, the assertion is verified.
\end{proof}

\begin{prop} \label{prop:Capq-Lq-weak}
Let $1 \leq q<\infty$, then the following statements are equivalent:
\begin{enumerate}
 \item
\begin{equation} \label{eq:Capq-implies-Lq-weak}
 \forall f \in \F \;\;\; D_1 \norm{f - M_\mu f}_{N(\mu),\infty} \leq
\norm{\abs{\nabla f}}_{L_q(\mu)} \;\; ,
\end{equation}
\item
\[
 Cap_q(t,1/2) \geq D_2 N^\wedge(t) \;\;\; \forall t \in [0,1/2] \;\; ,
\]
\end{enumerate}
and the best constants $D_1,D_2$ above satisfy $ D_1 \leq D_2 \leq 4 D_1$.
\end{prop}
\begin{proof}
$D_2 \geq D_1$ by Lemma \ref{lem:Lq-implies-capq}. To see the other
direction, note that as in Remark \ref{rem:cap-app0}, we may assume that $\int_{\set{f=t}}
\abs{\nabla f}^q d\mu = 0$ for all $t \in \Real$, and by replacing
$f$ with $f - M_\mu f$, that $M_\mu f = 0$. Note that if suffices to
show (\ref{eq:Capq-implies-Lq-weak}) with $D_1 = D_2$ for
non-negative functions for which $\mu\set{f=0} \geq 1/2$, since for
a general function as above, we can apply
(\ref{eq:Capq-implies-Lq-weak}) to $f_+ = f \chi_{f \geq 0}$ and to
$f_- = -f \chi_{f \leq 0}$, which yields:
\begin{multline*}
\norm{\abs{\nabla f}}_{L_q(\mu)} =
\brac{\int \abs{\nabla f_+}^q d\mu + \int \abs{\nabla f_-}^q d\mu}^{1/q} \geq
D_1 \brac{\norm{f_+}^q_{N(\mu),\infty} + \norm{f_-}^q_{N(\mu),\infty}}^{1/q} \\
\geq D_1 2^{1/q-1} \brac{
\norm{f_+}_{N(\mu),\infty} +
\norm{f_-}_{N(\mu),\infty}} \geq \frac{D_1}{4}
\norm{f}_{N(\mu),\infty} ~.
\end{multline*}

Given a non-negative function $f$ as above ($\mu\set{f=0} \geq 1/2$ hence $M_\mu f = 0$), and $t>0$, define $\Omega_t = \set{f \leq
t}$ and $f_t := f/t \wedge 1$. Then:
\begin{eqnarray*}
 \brac{\int_\Omega \abs{\nabla f}^q d\mu}^{1/q} &\geq& \brac{\int_{\Omega_t} \abs{\nabla f}^q
d\mu}^{1/q} \geq t \brac{\int_{\Omega} \abs{\nabla f_t}^q d\mu}^{1/q} \\
&\geq& t Cap_q(\mu\set{f_t \geq 1},1/2) \geq D_2 t N^\wedge(\mu\set{f \geq
t}) ~.
\end{eqnarray*}
Taking supremum on $t>0$, the assertion follows.
\end{proof}


\subsection{$q$-capacitary and strong Orlicz-Sobolev inequalities}

\begin{prop} \label{prop:Capq-Lq}
If $N(t)^{1/q}/t$ is non-decreasing on $\Real_+$ with $1 \leq q<\infty$, then the following statements are
equivalent:
\begin{enumerate}
 \item
\begin{equation} \label{eq:Capq-implies-Lq}
 \forall f \in \F \;\;\; D_1 \norm{f - M_\mu f}_{N(\mu)} \leq
\norm{\abs{\nabla f}}_{L_q(\mu)} \;\; ,
\end{equation}
\item
\[
 Cap_q(t,1/2) \geq D_2 N^\wedge(t) \;\;\; \forall t \in [0,1/2] \;\; ,
\]
\end{enumerate}
and the best constants $D_1,D_2$ above satisfy $ D_1 \leq D_2 \leq 4 D_1$.
\end{prop}
\begin{rem}
As already mentioned in Section \ref{sec:defs}, we call an inequality of the form
(\ref{eq:Capq-implies-Lq}) an Orlicz-Sobolev inequality (even though $N$ may not be convex).
\end{rem}
\begin{rem} \label{rem:improve}
One may show (see e.g. the proof of \cite[Theorem 1]{RobertoZegarlinski}) that when $N(t^{1/q})$ is
convex (so in particular $N(t)^{1/q}/t$ is non-decreasing),
Proposition \ref{prop:Capq-Lq} is equivalent to a theorem of Maz'ya
\cite[p. 112]{MazyaBook}, with a constant depending on $q$ which is better than the constant $4$ above. In particular, when $q=1$, Maz'ya showed that the optimal constant is actually $1$, so that $D_1 = D_2$. The latter conclusion was also independently derived by Federer--Fleming \cite{FedererFleming}.
\end{rem}

\begin{proof}
$D_2 \geq D_1$ by (\ref{eq:weak-norm}) and Lemma
\ref{lem:Lq-implies-capq}. To see the other direction, we assume
again (as in Remark \ref{rem:cap-app0}) that $\int_{\set{f=t}}
\abs{\nabla f}^q d\mu = 0$ for all $t \in \Real$, and by replacing
$f$ with $f - M_\mu f$, that $M_\mu f = 0$. Again, if suffices to
show (\ref{eq:Capq-implies-Lq}) for non-negative functions for which
$\mu\set{f=0} \geq 1/2$, but now we do not lose in the constant.
Indeed, for a general function as above, we can apply
(\ref{eq:Capq-implies-Lq}) to $f_+ = f \chi_{f \geq 0}$ and to $f_-
= -f \chi_{f \leq 0}$, which yields:
\[
\norm{\abs{\nabla f}}^q_{L_q(\mu)} =
\int \abs{\nabla f_+}^q d\mu + \int \abs{\nabla f_-}^q d\mu \geq
D_1^q \brac{\norm{f_+}^q_{N(\mu)} + \norm{f_-}^q_{N(\mu)}} \geq
 D_1^q \norm{f}^q_{N(\mu)}.
\]
The last inequality follows from the fact that 
$N^{1/q}(t)/t$ is non-decreasing, so denoting $v_{\pm} = \norm{f_{\pm}}_{N(\mu)}$,
we indeed verify that:
\begin{multline*}
\int N\brac{\frac{f_+ + f_-}{(v^q_+ + v^q_-)^{1/q}}} d\mu =
 \int N\brac{\frac{f_+}{v_+} \frac{v_+}{(v^q_+ + v^q_-)^{1/q}}} d\mu +
\int N\brac{\frac{f_-}{v_-} \frac{v_-}{(v^q_+ + v^q_-)^{1/q}}} d\mu \\
\leq
\frac{v_+^q}{v^q_+ + v^q_-}  \int  N\brac{\frac{f_+}{v_+}}d\mu +
\frac{v_-^q}{v^q_+ + v^q_-}  \int  N\brac{\frac{f_-}{v_-}}d\mu \leq 1 ~.
\end{multline*}

We will first assume that $f$ is bounded. Given a bounded
non-negative function $f$ as above ($M_\mu f = 0$ and $\mu\set{f=0}
\geq 1/2$), we may assume by homogeneity that $\norm{f}_{L_\infty} =
1$. For $i\geq 1$, denote $\Omega_i = \set{1/2^{i} \leq f \leq
1/2^{i-1}}$, $m_i = \mu(\Omega_i)$, $f_i = 2^{i}(f - 1/2^{i}) \vee 0
\wedge 1$ and set $m_0=0$. Also denote $J := N^\wedge$. Now:
\begin{multline*}
\norm{\abs{\nabla f}}_{L_q(\mu)}^q = \sum_{i=1}^\infty
\int_{\Omega_i} \abs{\nabla f}^q d\mu \geq
\sum_{i=1}^\infty \frac{1}{2^{qi}} \int_{\Omega} \abs{\nabla f_i}^q d\mu \\
\geq \sum_{i=1}^\infty
\frac{1}{2^{q i}} Cap^q_q(\mu\set{f \geq 1/2^{i-1}},1/2) \geq D_2^q \sum_{i=2}^\infty
\frac{J^q(m_{i-1})}{2^{q i}} = \frac{D_2^q}{4^q} V^q,
\end{multline*}
where:
\[
V := \brac{\sum_{i=1}^\infty \frac{J^q(m_{i})}{2^{q (i-1)}}}^{1/q} ~.
\]
It remains to show that $\norm{f}_{N(\mu)} \leq V$. Indeed:
\[
 \int_\Omega N\brac{\frac{f}{V}} d\mu \leq \sum_{i=1}^\infty m_i N\brac{\frac{1}{2^{i-1}
V}}
= \sum_{i=1}^\infty \frac{J^{-1}(J(m_i))}{J^{-1}(2^{i-1} V)} \leq \sum_{i=1}^\infty
\frac{J^q(m_i)}{2^{q(i-1)} V^q} = 1,
\]
where in the last inequality we have used the fact that $N(t)^{1/q}/t$ is non-decreasing, hence
$(J^{-1})^{1/q}(t)/t$ is non-decreasing, and therefore:
\[
\frac{ J^{-1}(x) } {J^{-1}(y) } \leq \brac{\frac{x}{y}}^q,
\]
whenever $x/y \leq 1$, which is indeed the case for us.

For a non-bounded $f\in \F$ with $\mu\set{f=0} \geq 1/2$, we may
define $f_m = f \wedge b_m$ so that $\mu\set{f
> b_m} \leq 1/m$ and (just for safety) $\mu\set{f=b_m} = 0$. It then follows
by what was proved for bounded functions that:
\[
\norm{\abs{\nabla f}}_{L_q(\mu)} \geq \lim_{m \rightarrow \infty}
\norm{\abs{\nabla f_m}}_{L_q(\mu)} \geq D_1 \lim_{m \rightarrow
\infty} \norm{f_m}_{N(\mu)} = D_1 Z ~,
\]
where all limits exist since they are non-decreasing. To conclude,
$Z \geq \norm{f}_{N(\mu)}$, since $N$ is continuous, so by
the Monotone Convergence Theorem:
\[
\int N(f/Z) d\mu = \int \lim_{m \rightarrow \infty}
N(f_m/Z) d\mu = \lim_{m \rightarrow \infty} \int
N(f_m/Z) d\mu \leq 1 ~.
\]
\end{proof}


\subsection{Passing between $q$-capacitary inequalities}

The case $q_0=1$ in the following proposition is due to Maz'ya
\cite[p. 105]{MazyaBook}. Following \cite{EMilmanRoleOfConvexityInFunctionalInqs}, we provide a proof which generalizes to the case of an arbitrary metric probability space and $q_0 > 1$. We denote the conjugate exponent to
$q \in [1,\infty]$ by $q^* = q/(q-1)$.

\begin{prop} \label{prop:increase-cap-q}
Let $1 \leq q_0 \leq q < \infty$ and set $p_0 = q_0^*, p = q^*$. Then for all $0<a<b<1$:
\[
 \frac{1}{Cap_q(a,b)} \leq \gamma_{p,p_0} \brac{\int_a^b
\frac{ds}{(s-a)^{p/p_0} Cap_{q_0}^p(s,b)}}^{1/p} ~,
\]
where:
\begin{equation} \label{eq:gamma}
\gamma_{p,p_0} :=  \frac{(\frac{p_0}{p} - 1)^{1/p_0}}{(1-\frac{p}{p_0})^{1/p}} ~.
\end{equation}
\end{prop}

\begin{proof}
Let $0<a<b<1$ be given, and let $\Phi : \Omega \rightarrow [0,1]$ be
a function in $\F$ such that $a':=\mu\set{\Phi=1} \geq a$ and
$1-b':=\mu\set{\Phi=0} \geq 1-b$. As usual (see Remark
\ref{rem:cap-app0}), by approximating $\Phi$, we may assume that
$\int_{\set{\Phi=t}} \abs{\nabla \Phi}^q d\mu = 0$ for all $t \in
(0,1)$. Let $C := \set{ t \in (0,1) ; \mu\set{\Phi = t} > 0 }$
denote the discrete set of atoms of $\Phi$ under $\mu$, set $\Gamma
:= \set{ f \in C }$ and denote $\gamma = \mu(\Gamma)$.

We now choose $t_0=0 < t_1 < t_2 < \ldots < 1$, so that denoting for
$i \geq 1$, $\Omega_i = \set{t_{i-1} \leq \Phi \leq t_i}$, and
setting $m_i = \mu(\Omega_i \setminus \Gamma)$, we have $m_i =
(b'-a'-\gamma) \alpha^{i-1}(1-\alpha)$, where $0 \leq \alpha \leq 1$
will be chosen later. Denote in addition $\Phi_i =
\brac{\frac{\Phi-t_{i-1}}{t_i - t_{i-1}} \vee 0} \wedge 1$,
$N_i = \sum_{j>i} m_j$. Applying H\"{o}lder's
inequality twice, we estimate:
\begin{eqnarray*}
\brac{\int_\Omega \abs{\nabla \Phi}^q d\mu}^{1/q} &=&
\brac{\sum_{i=1}^\infty \int_{\Omega_i \setminus \Gamma} \abs{\nabla
\Phi}^q d\mu}^{1/q} \geq \brac{\sum_{i=1}^\infty
m_i^{1-\frac{q}{q_0}} \brac{\int_{\Omega_i \setminus \Gamma}
\abs{\nabla \Phi}^{q_0}
d\mu}^{q/q_0}}^{1/q} \\
&\geq& \brac{\sum_{i=1}^\infty m_i^{1-\frac{q}{q_0}} (t_i-t_{i-1})^q
\brac{\int_{\Omega} \abs{\nabla \Phi_i}^{q_0} d\mu}^{q/q_0}}^{1/q} \\
&\geq& \brac{\sum_{i=1}^\infty m_i^{1-\frac{q}{q_0}}
(t_i-t_{i-1})^q Cap_{q_0}^q(\mu\set{\Phi_i = 1},1-\mu\set{\Phi_i = 0})}^{1/q} \\
&\geq& \sum_{i=1}^\infty (t_i-t_{i-1}) \brac{\sum_{i=1}^\infty
\frac{m_i^{1-p/p_0}}{Cap_{q_0}^p(\mu\set{\Phi \geq t_i},b)}}^{-1/p}
~.
\end{eqnarray*}
Since $\mu\set{\Phi \geq t_i} \geq a' + N_i$ and $Cap_{q_0}(s,b)$ is
non-decreasing in $s$, we continue to estimate as follows:
\begin{eqnarray*}
& & \brac{\frac{1}{\int_\Omega \abs{\nabla \Phi}^q d\mu}}^{p/q} \leq \sum_{i=1}^\infty
\frac{m_i^{1-p/p_0}}{Cap_{q_0}^p(\mu\set{\Phi \geq t_i},b)} \\
&\leq & \sum_{i=1}^\infty \frac{m_i^{1-p/p_0}}{m_{i+1}}
\int_{a'+N_{i+1}}^{a'+N_i} \frac{ds}{Cap_{q_0}^p(s,b)} \leq
\sum_{i=1}^\infty \frac{1}{\alpha m_i^{p/p_0}}
\int_{a'+N_{i+1}}^{a'+N_i} 
\frac{ds}{Cap_{q_0}^p(s,b)} \\
&\leq & \frac{1}{\alpha} \brac{\frac{\alpha}{1-\alpha}}^{p/p_0}  \sum_{i=1}^\infty
\int_{a'+N_{i+1}}^{a'+N_i} 
\frac{ds}{(s-a')^{p/p_0} Cap_{q_0}^p(s,b)} \\
&\leq& \frac{1}{\alpha} \brac{\frac{\alpha}{1-\alpha}}^{p/p_0}  \int_a^b \frac{ds}{(s-a)^{p/p_0}
Cap_{q_0}^p(s,b)} ~,
\end{eqnarray*}
where we have used that $m_{i+1} = \alpha m_i$, $m_i =
\frac{1-\alpha}{\alpha} N_i$, and in the last inequality the fact
that $Cap_{q_0}(s,b)$ is non-decreasing in $s$. The assertion now
follows by taking supremum on all $\Phi$ as above, and choosing the
optimal $\alpha = 1 - p/p_0$.
\end{proof}


\subsection{Combining everything}

Combining all of the ingredients in this section, we see how to pass from isoperimetric inequalities to functional inequalities, simply by following the general diagram:
\begin{eqnarray*}
\text{Isoperimetric inequality} & \Leftrightarrow_{\text{Corollary \ref{cor:cap1}}} & \text{$1$-capacitary inequality} \\
& & \Downarrow \text{\begin{small}Proposition \ref{prop:increase-cap-q}\end{small}} \\
\begin{array}[t]{c} \text{$(N,q)$ Orlicz-Sobolev inequality} \\ \text{with $N(t)^{1/q}/t$ non-decreasing} \end{array} & \Leftrightarrow_{\text{Proposition \ref{prop:Capq-Lq}}} & \text{$q$-capacitary inequality}
\end{eqnarray*}

In particular, it is an exercise to follow this diagram and obtain the previously mentioned inequalities of Subsection \ref{subsec:known-connections}: $D_{Poin} \geq c D_{Lin}$, $D_{LS_2} \geq c D_{Gau}$ and $D_{LS_q} \geq c D_{Exp_q}$ for $q \in [1,2]$, for some universal constant $c>0$. In the first case, the optimal constant $c=1/2$ may also be obtained by improving the constant in Proposition \ref{prop:Capq-Lq} as in Remark \ref{rem:improve} (but of course easier ways are known to obtain this optimal constant, see e.g. \cite{EMilman-RoleOfConvexity}). More generally, the following statement may easily be obtained (see \cite{EMilmanRoleOfConvexityInFunctionalInqs} for more details and useful special cases):

\begin{thm} \label{thm:I-implies-Orlicz}
Let $1 \leq q < \infty$, and set $p = q^*$. Let $N \in \J$, so that
$N(t)^{1/q}/t$ is non-decreasing. Then:
\begin{equation} \label{eq:I-implies-Orlicz-assumption}
 \tilde{I}(t) \geq D t^{1-1/q} N^\wedge(t) \;\; \forall t \in [0,1/2]
\end{equation}
implies:
\begin{equation} \label{eq:I-implies-Orlicz-conclusion}
\forall f \in \F \;\; B_{N,q} D \norm{f - M_\mu f}_{N(\mu)} \leq \norm{ \abs{\nabla f}}_{L_q(\mu)}
~,
\end{equation}
where:
\begin{equation} \label{eq:B-estimate}
 B_{N,q} := \frac{1}{4} \inf_{0<t<1/2} \frac{1}{\brac{\int_t^{1/2}
\frac{N^\wedge(t)^p ds}{s N^\wedge(s)^p  }}^{1/p}} ~.
\end{equation}
\end{thm}


\section{The Converse Statement} \label{sec:new}

Our goal in this section will be to prove the following converse to Theorem \ref{thm:I-implies-Orlicz}:

\begin{thm} \label{thm:Orlicz-implies-I}
Let $1 < q \leq \infty$ and let $N \in \J$ denote a Young function so that $N(t)^{1/q}/t$ is non-decreasing. 
Then under our $\kappa$-semi-convexity
assumptions, the statement:
\begin{equation} \label{eq:Orlicz-inq}
 \forall f \in \F \;\; D \norm{f - M_\mu f}_{N(\mu)} \leq \norm{ \abs{\nabla f} }_{L_q(\mu)}
\end{equation}
implies:
\begin{equation} \label{eq:Orlicz-implies-I-conclusion}
 \tilde{I}(t) \geq C_{N,q} \min\brac{D,\frac{D^{r}}{\kappa^{\frac{r-1}{2}}}} t^{1-1/q} N^\wedge(t) \;\; \forall t \in [0,1/2] ~,
\end{equation}
where $r = \max(q,2)$ and $C_{N,q} > 0$ depends solely on $N$ and $q$.
\end{thm}

\begin{rem}
The assumption that $N$ is a convex function is not essential for this result, since it is possible to approximate $N$ appropriately in the large using a convex function as in Subsection \ref{subsubsec:capacity-approach} below. We refer to \cite[Theorem 4.5]{EMilmanRoleOfConvexityInFunctionalInqs} for more details. 
\end{rem}

\begin{rem}
Clearly, using the results of Section \ref{sec:caps}, Theorem \ref{thm:Orlicz-implies-I} may be reformulated as a converse to Maz'ya's inequality relating $Cap_q$ and $Cap_1$ (Proposition \ref{prop:increase-cap-q}), under our $\kappa$-semi-convexity assumptions. For the case of our convexity assumptions ($\kappa=0$), this has been explicitly written out in \cite[Theorem 5.1]{EMilmanRoleOfConvexityInFunctionalInqs}.  
\end{rem}

The case $\kappa = 0$ of Theorem \ref{thm:Orlicz-implies-I} was proved in \cite{EMilmanRoleOfConvexityInFunctionalInqs} using the semi-group approach developped by Bakry--Ledoux \cite{BakryLedoux} and Ledoux \cite{LedouxBusersTheorem,LedouxSpectralGapAndGeometry}. Let us now recall this framework. 

Given a smooth complete connected Riemannian manifold $\Omega = (M,g)$ equipped with a
probability measure $\mu$ with density $d\mu = \exp(-\psi) dvol_M$, $\psi \in C^2(M,\Real)$, we
define the associated Laplacian $\Delta_{(\Omega,\mu)}$ by:
\begin{equation} \label{eq:Laplacian-def}
 \Delta_{(\Omega,\mu)} := \Delta_{\Omega} - \nabla \psi \cdot \nabla,
\end{equation}
where $\Delta_{\Omega}$ is the usual Laplace-Beltrami operator on $\Omega$. $\Delta_{(\Omega,\mu)}$
acts on $\B(\Omega)$, the space of bounded smooth real-valued functions on $\Omega$.
Let $(P_t)_{t \geq 0}$ denote the semi-group associated to the diffusion process
with infinitesimal generator $\Delta_{(\Omega,\mu)}$
(cf. \cite{DaviesSemiGroupBook,LedouxLectureNotesOnDiffusion}), for which $\mu$ is its stationary measure.  
It is characterized by the following system of second order differential equations:
\begin{equation} \label{eq:characterize-Delta}
 \frac{d}{dt} P_t(f) = \Delta_{(\Omega,\mu)} (P_t(f)) \;\;\;\; P_0(f) = f \;\;\; \forall f \in
\B(\Omega)~.
\end{equation}
For each $t \geq 0$, $P_t : \B(\Omega) \rightarrow \B(\Omega)$ is a bounded linear operator in the $L_\infty$ norm, and
its action naturally extends to the entire $L_p(\mu)$ spaces ($p\geq 1$). We collect several
elementary properties of these operators:
\begin{itemize}
\item
$P_t 1 = 1$.
 \item
$f \geq 0 \Rightarrow P_t f \geq 0$.
\item
$\int (P_t f) g d\mu = \int f (P_t g) d\mu$.
\item
$N(\abs{P_t(f)}) \leq P_t(N(\abs{f}))$ for any Young function $N$.
\item
$P_t \circ P_s = P_{t+s}$.
\end{itemize}

The following crucial dimension-free reverse Poincar\'e inequality was shown by Bakry--Ledoux in
\cite[Lemma 4.2]{BakryLedoux}, extending Ledoux's approach \cite{LedouxBusersTheorem} for proving
Buser's Theorem (see also \cite[Lemma 2.4]{BakryLedoux}, \cite[Lemma 5.1]{LedouxSpectralGapAndGeometry}):

\begin{lem}[Bakry--Ledoux] \label{lem:Ledoux}
Assume that the following Bakry--\'{E}mery Curvature-Dimension condition holds on $\Omega$:
\begin{equation} \label{eq:Ledoux-condition}
Ric_g + Hess_g \psi \geq -\kappa g ~, \kappa \geq 0 ~.
\end{equation}
Then for any $t\geq 0$ and $f \in \B(\Omega)$, we have:
\[
 K(\kappa,t) \abs{\nabla P_t f}^2 \leq P_t(f^2) - (P_t f)^2
\]
pointwise, where:
\[
 K(\kappa,t) := \frac{1 - \exp(-2\kappa t)}{\kappa} \;\;\; \brac{= 2t \;\; \text{ if } \kappa =0} ~.
\]
\end{lem}
\begin{proof}[Sketch of Proof following \cite{LedouxSpectralGapAndGeometry}]
Note that $\Delta_{(\Omega,\mu)}(g^2) - 2 g \Delta_{(\Omega,\mu)}(g) = 2 |\nabla g|^2$ for any $g \in \B(\Omega)$. Consequently, (\ref{eq:characterize-Delta}) implies: 
\[
 P_t(f^2) - (P_t f)^2 = \int_0^t \frac{d}{ds} P_s((P_{t-s}f)^2) ds = 2 \int_0^t P_s(| \nabla P_{t-s} f |^2) ds ~.
\]
The main observation is that under the Bakry-\'Emery condition (\ref{eq:Ledoux-condition}), the function $\exp(2\kappa s) P_s(| \nabla P_{t-s} f |^2)$ is non-decreasing, as verified by direct differentiation and use of the Bochner formula. Therefore:
\[
  P_t(f^2) - (P_t f)^2 \geq |\nabla P_t f|^2 \; 2 \int_0^t \exp(-2\kappa s) ds ~,
\]
which concludes the proof.
\end{proof}

In fact, the proof of this lemma is very general and extends to the abstract
framework of diffusion generators, as developed by Bakry and
\'{E}mery in their celebrated paper \cite{BakryEmery}. In the
Riemannian setting, it is known \cite{QianGradientEstimateWithBoundary} (see also \cite{HsuGradientEstimateWithBoundary,WangGradientEstimateWithBoundary}) that the gradient estimate
of Lemma \ref{lem:Ledoux} remains valid when the support of $\mu$ is the closure of a
locally convex domain (connected open set) $\Omega \subset (M,g)$ with $C^2$ boundary, and
$d\mu|_\Omega = \exp(-\psi) dvol_M|_\Omega$, $\psi \in C^2(\overline{\Omega},\Real)$. A domain $\Omega$ with $C^2$ boundary is called locally convex if the second fundamental form on $\partial \Omega$ is positive semi-definite (with respect to the normal field pointing inward). In this case, $\Delta_{\Omega}$ in (\ref{eq:Laplacian-def}) denotes the Neumann Laplacian on $\overline{\Omega}$, $\B(\Omega)$ denotes the space of bounded smooth real-valued functions on $\overline{\Omega}$ satisfying
Neumann's boundary conditions on $\partial \Omega$, and Lemma \ref{lem:Ledoux} remains valid.

Under these assumptions, Lemma \ref{lem:Ledoux} clearly implies that:
\begin{equation} \label{eq:L_q-bound}
 \forall q \in [2,\infty] \;\;\; \forall f \in \B(\Omega) \;\;\; \norm{\abs{\nabla P_t f} }_{L_q(\mu)} \leq \frac{1}{\sqrt{K(\kappa,t)}} \norm{f}_{L_q(\mu)} ~,
\end{equation}
and using $q = \infty$, Ledoux easily deduces the following dual statement
(cf. \cite[(5.5)]{LedouxSpectralGapAndGeometry}):
\begin{lem}[Ledoux] \label{lem:Ledoux-dual}
Under the same assumptions as Lemma \ref{lem:Ledoux}, we have:
\begin{equation} \label{eq:L_1-dual-bound}
 \forall f \in \B(\Omega) \;\;\; \norm{ f - P_t f}_{L_1(\mu)} \leq \int_0^t \frac{ds}{\sqrt{K(\kappa,s)}} \norm{ \abs{\nabla f} }_{L_1(\mu)} ~.
\end{equation}
\end{lem}
\noindent
To use (\ref{eq:L_q-bound}) and (\ref{eq:L_1-dual-bound}), it will be convenient to note the following rough estimates:
\begin{equation} \label{eq:rough}
  t \in [0,1/(2\kappa)] \;\;\; \Rightarrow \;\; K(\kappa,t) \geq t \;\;\;,\;\;\; \int_0^t \frac{ds}{\sqrt{K(\kappa,s)}} \leq 2 \sqrt{t} ~.
\end{equation}

It will also be useful to introduce the following:

\begin{dfn*}
We denote by $N(\mu)^*$ the dual norm to $N(\mu)$, given by:
\[
 \norm{f}_{N(\mu)^*} := \sup \set{ \int f g d\mu ; \norm{g}_{N(\mu)} \leq 1}.
\]
\end{dfn*}

\noindent
It is elementary to calculate the $N(\mu)^*$-norm of characteristic functions (cf. \cite[p. 111]{MazyaBook}):

\begin{lem} \label{lem:Mazya-duality}
Let $N$ denote a Young function. Then for any Borel set $A$ with $\mu(A)>0$:
\[
\norm{\chi_A}_{N(\mu)^*} = \mu(A) N^{-1}\brac{\frac{1}{\mu(A)}} = \frac{\mu(A)}{N^\wedge(\mu(A))} ~.
\]
\end{lem}

\subsection{Case of $q \geq 2$}

To handle the $\kappa > 0$ case, we will need the following new estimate, which may be of independent interest:

\begin{prop}\label{prop:new}
Assume that Bakry--\'{E}mery Curvature-Dimension condition (\ref{eq:Ledoux-condition}) holds on $\Omega$, 
and that the following $(N,q)$ Orlicz-Sobolev inequality is satisfied for $N \in \J$ and $q \geq 2$:
\begin{equation} \label{eq:Orlicz-inq-E}
 \forall f \in \F \;\;\;\; D \norm{f - E_\mu f}_{N(\mu)} \leq \norm{ \abs{\nabla f} }_{L_q(\mu)} ~.
\end{equation}
Then for any $f \in L_\infty(\Omega)$ with $\int f d\mu = 0$, we have for all $t \geq 0$:
\[
\int |P_t f|^2 d\mu \leq \int f^2 d\mu \brac{1 + (q-1) \frac{2 D^q}{\norm{f}^q_{N(\mu)^*}\norm{f}^{q-2}_{L_\infty}} \brac{\int f^2 d\mu}^{q-1} \int_0^t K(\kappa,s)^{\frac{q-2}{2}} ds  }^{-\frac{1}{q-1}} ~.
\]
\end{prop}
\begin{proof}
Denote $u(t):= \int |P_t f|^2 d\mu = \int f P_{2t}f d\mu$ by self-adjointness and the semi-group property. 
By (\ref{eq:characterize-Delta}) and integration by parts, we have:
\begin{equation} \label{eq:parts}
 u'(t) = 2 \int P_t f \Delta_{(\Omega,\mu)}(P_t f) d\mu = -2 \int | \nabla P_t f |^2 d\mu ~.
\end{equation} 
Note that $u(t)$ is decreasing. We now use the following estimate:
\begin{eqnarray*}
u(t)^q &\leq & u(t/2)^q = 
 \brac{\int f P_t f d\mu}^q \leq \norm{f}^q_{N(\mu)^*} \norm{P_{t} f}^q_{N(\mu)}  \\
&\leq&  \frac{\norm{f}^q_{N(\mu)^*} }{D^q} \int |\nabla P_{t} f|^q d\mu 
\leq  \frac{\norm{f}^q_{N(\mu)^*}}{D^q}  \int |\nabla P_{t} f|^2 d\mu \norm{|\nabla P_{t} f|}_{L_\infty}^{q-2} ~.
\end{eqnarray*}
By (\ref{eq:L_q-bound}) and (\ref{eq:parts}), we obtain:
\[
 u(t)^q \leq - \frac{\norm{f}^q_{N(\mu)^*}}{2 D^q} \frac{\norm{f}^{q-2}_{L_\infty}}{K(\kappa,t)^{\frac{q-2}{2}}} \; u'(t) ~. 
\]
Denoting $v(t) = u(t)^{-q+1}$, we see that this boils down to:
\[
 v'(t) \geq (q-1) \frac{2 D^q}{\norm{f}^q_{N(\mu)^*}\norm{f}^{q-2}_{L_\infty}} K(\kappa,t)^{\frac{q-2}{2}} ~.
\]
Integrating in $t$, the desired conclusion follows.
\end{proof}

\begin{rem}
When $N(x) = x^2$, or more generally, when $N(x^{1/2})$ is convex, it is possible to obtain a better dependence on $q$ in Proposition \ref{prop:new}, which, as $q \rightarrow 2$, would recover the exponential rate of convergence of $\int |P_t f|^2 d\mu$ to $0$, as dictated by the spectral theorem. Unfortunately, it seem that this would not yield the correct dependence in $N$ in the assertion of Theorem \ref{thm:Orlicz-implies-I}.  
\end{rem}

\begin{proof}[Proof of Theorem \ref{thm:Orlicz-implies-I} for $q \geq 2$]

We will prove the theorem under the \emph{smooth} $\kappa$-semi-convexity assumptions of this section.
The general case follows by an approximation argument which was derived in \cite[Section 6]{EMilmanRoleOfConvexityInFunctionalInqs} for the 
case $\kappa = 0$, but holds equally true for any $\kappa \geq 0$. 

Since $N$ is a Young function, we may invoke Lemma \ref{lem:E-M} and replace $M_\mu f$ in (\ref{eq:Orlicz-inq})
by $E_\mu f$ as in (\ref{eq:Orlicz-inq-E}), at the expense of an additional universal constant in
the final conclusion.

Let $A$ denote an arbitrary Borel set in $\Omega$ so that $\mu^+(A) < \infty$, and let $\chi_{A,\eps,\delta}(x) :=
(1 - \frac{1}{\eps} d(x,A^d_{\delta})) \vee 0$ be a continuous approximation in $\Omega$ to the
characteristic function $\chi_A$ of $A$ (as usual $d$ denotes the induced geodesic distance). Our assumptions imply that:
\[
 \frac{\mu(A^d_{\eps+2\delta}) - \mu(A)}{\eps} \geq \int \abs{\nabla{\chi_{A,\eps,\delta}}} d\mu ~.
\]
Applying Lemma \ref{lem:Ledoux-dual} to functions in $\B(\Omega)$ which approximate $\chi_{A,\eps,\delta}$
(in say $W^{1,1}(\Omega,\mu)$) and passing to the limit
inferior as $\eps,\delta \rightarrow 0$ so that $\delta / \eps \rightarrow 0$, it follows that:
\[
\int_0^t \frac{ds}{\sqrt{K(\kappa,s)}} \mu^+(A) \geq \int \abs{\chi_A - P_t \chi_A} d\mu 
\]
(note that the assumption $\mu^+(A) < \infty$ guarrantees that $\mu(\overline{A} \setminus A) =0$, so $\chi_{A,\eps,\delta}$ tends to $\chi_A$ in $L_1(\mu)$).  
We start by rewriting the right hand side above as:
\begin{eqnarray*}
& & \!\!\!\!\!\!\!\!\!\!\!\!\!\!\!\!\!\! \int_A (1 - P_t \chi_A) d\mu + \int_{\Omega \setminus A}
P_t \chi_A d\mu = 2\brac{\mu(A) - \int_A
P_t \chi_A d\mu} \\
& = & 2 \brac{\mu(A)(1-\mu(A)) - \int_\Omega (P_t \chi_A - \mu(A))(\chi_A - \mu(A)) d\mu} \\
& = & 2 \brac{\int_\Omega |\chi_A - \mu(A)|^2 d\mu - \int_\Omega |P_{t/2}(\chi_A - \mu(A))|^2 d\mu}  ~.
\end{eqnarray*}
Denoting $f = \chi_A - \mu(A)$ and using Proposition \ref{prop:new} to estimate the right-most expression, we 
obtain after using the estimates in (\ref{eq:rough}), that for $t \leq 1/(2\kappa)$:
\begin{equation} \label{eq:optimize-on}
2 \sqrt{t} \mu^+(A) \geq 2 \mu(A) (1-\mu(A)) \brac{1 - \brac{1 + (q-1) M_t}^{-\frac{1}{q-1}}} ~,
\end{equation}
where:
\[
M_t := \frac{2 D^q}{\norm{f}^q_{N(\mu)^*}\norm{f}^{q-2}_{L_\infty}} \brac{\int f^2 d\mu}^{q-1} \frac{2}{q} \brac{\frac{t}{2}}^{\frac{q}{2}}  ~.
\]

To estimate $M_t$, we employ Lemma \ref{lem:Mazya-duality}:
\begin{eqnarray*}
& & \!\!\!\!\!\!\!\!\!\!\!\!\!\! \norm{\chi_A - \mu(A)}_{N(\mu)^*} \leq (1-\mu(A)) \norm{\chi_A
}_{N(\mu)^*} + \mu(A)\norm{\chi_{\Omega\setminus A} }_{N(\mu)^*} \\
&=& \mu(A)(1-\mu(A))  \brac{\frac{1}{N^{\wedge}(\mu(A))} + \frac{1}{N^{\wedge}(1-\mu(A))}}  \leq
2 \frac{\mu(A)(1-\mu(A)) } {N^\wedge(\min(\mu(A),1-\mu(A)))} ~,
\end{eqnarray*}
and using that $\int f^2 d\mu = \mu(A)(1-\mu(A))$, we conclude that:
\begin{equation} \label{eq:new-M_t-estimate}
 M_t \geq L_t := E \; \frac{2}{q} \brac{\frac{t}{2}}^{\frac{q}{2}} ~,~ E := \frac{D^q N^\wedge(\min(\mu(A),1-\mu(A)))^q}{2^{q-1} \mu(A) (1-\mu(A))} 
\end{equation}

It remains to optimize on $t$ in (\ref{eq:optimize-on}). Denote:
\[
 t_0 := 4 \brac{\frac{q}{2 E}}^{2/q} = 4 \brac{\frac{q}{2}}^{2/q} \frac{2^{2/q^*} (\mu(A)(1-\mu(A)))^{2/q}}{D^2 N^\wedge(\min(\mu(A),1-\mu(A)))^2} ~.
\]
\begin{enumerate}
\item 
If $t_0 \leq 1/(2\kappa)$, we see from (\ref{eq:new-M_t-estimate}) that $L_{t_0} = 2^{q/2}$,
and we immediately obtain from (\ref{eq:optimize-on}) that:
\[
 \mu^+(A) \geq c_1 \frac{\mu(A) (1-\mu(A))}{\sqrt{t_0}} \geq c_2 D \min(\mu(A),1-\mu(A))^{1-1/q} N^\wedge(\min(\mu(A),1-\mu(A))) ~, 
\]
where $c_1,c_2>0$ are some numeric constants. 
\item
If $t_0 > 1/(2\kappa)$, we evaluate (\ref{eq:new-M_t-estimate}) and (\ref{eq:optimize-on}) at time $t_1 = 1/(2\kappa)$, for which $L_{t_1} < 2^{q/2}$. Therefore $(1+(q-1) L_{t_1})^{\frac{1}{q-1}} \geq 1+ c_3 2^{-q/2} L_{t_1}$ (recall that $q \geq 2$), and hence by (\ref{eq:optimize-on}):
\[
 \mu^+(A) \geq c_4 \sqrt{\kappa} \mu(A) (1-\mu(A)) 2^{-q/2} L_{t_1} ~, 
\]
where $c_3,c_4 > 0$ are numeric constants. Plugging in $L_{t_1}$, we obtain:
\[
 \mu^+(A) \geq \frac{c_5}{q 4^q} \frac{D^q}{\kappa^{\frac{q-1}{2}}} N^\wedge(\min(\mu(A),1-\mu(A)))^q ~.
\]
Using that $N(t)^{1/q}/t$ is non-decreasing, which is equivalent to $N^\wedge(t)/t^{1/q}$ being non-increasing, we conclude that:
\begin{equation} \label{eq:degrade-1}
 \mu^+(A) \geq \frac{c_6 N^\wedge(1/2)^{q-1}}{q 4^q} \frac{D^q}{\kappa^{\frac{q-1}{2}}} \min(\mu(A),1-\mu(A))^{1-1/q} N^\wedge(\min(\mu(A),1-\mu(A))) ~,
\end{equation}
for some numeric constants $c_5,c_6>0$.
\end{enumerate}
Combining both cases, the assertion follows in the case $q \geq 2$.
\end{proof}

We conclude the study of the case $q \geq 2$ by mentioning that, even though our estimates in (\ref{eq:degrade-1}) degrade as $q \rightarrow \infty$, it is also possible to study the limiting case $q = \infty$. In this case, the functional inequality (\ref{eq:Orlicz-inq}) corresponds to an integrability property of \emph{Lipschitz} functions $f$ with $M_\mu(f) =0$, or equivalently, to the \emph{concentration} of the measure $\mu$, in terms of the decay of $\mu(\Omega \setminus A^d_t)$ as a function of $t$ for sets $A$ with $\mu(A) \geq 1/2$. Using techniques from Riemannian Geometry, it is still possible to deduce in this case an appropriate isoperimetric inequality under our $\kappa$-semi-convexity assumptions (and an appropriate necessary assumption on the concentration), see \cite{EMilmanGeometricApproachPartI}. 

\subsection{Case of $1 < q \leq 2$}

We will present two proofs of Theorem \ref{thm:Orlicz-implies-I} for this case, each having its own advantages and drawbacks.
The first runs along the same lines as in the previous subsection, and is new even in the $\kappa=0$ case. It has the advantage of working for arbitrary $N \in \J$ satisfying the assumptions of Theorem \ref{thm:Orlicz-implies-I}, but with this approach the estimates on $C_{N,q}$ degrade as $q \rightarrow 1$. This does not necessarily happen with the second proof, which is based on the idea in \cite{EMilmanRoleOfConvexityInFunctionalInqs} of reducing the claim to the $q=2$ case using Proposition \ref{prop:increase-cap-q}. Nevertheless, some further conditions on $N$ will need to be imposed for this approach to work.

\subsubsection{Semi-Group Approach}

We will need an analogue of Proposition \ref{prop:new}, which again may be of independent interest:

\begin{prop}\label{prop:new2}
Assume that the following $(N,q)$ Orlicz-Sobolev inequality is satisfied for $N \in \J$ and $1 < q \leq 2$:
\begin{equation} \label{eq:Orlicz-inq-E2}
 \forall f \in \F \;\;\;\; D \norm{f - E_\mu f}_{N(\mu)} \leq \norm{ \abs{\nabla f} }_{L_q(\mu)} ~.
\end{equation}
Then for any $f \in L_\infty(\Omega)$ with $\int f d\mu = 0$, we have for all $t \geq 0$:
\[
\int |P_t f|^q d\mu \leq \int f^q d\mu \brac{1 + \frac{2 D^2}{\norm{f}^{\frac{2(2-q)}{q-1}}_{L_1(\mu)} \norm{f}^2_{N(\mu)^*}} \brac{\int f^q d\mu}^{\frac{2}{q(q-1)}} t }^{-\frac{q(q-1)}{2}} ~.
\]
\end{prop}

Note that the case $q=2$ is identical to the one in Proposition \ref{prop:new}.

\begin{proof}
Denote $u(t) = \int |P_t f|^q d\mu$, and observe that:
\[
 u'(t) = q \int |P_t f|^{q-1} sign(P_t f) \Delta_{\Omega,\mu}(P_t f) d\mu = - q(q-1) \int |P_t f|^{q-2} |\nabla P_t f|^2 d\mu ~. 
\]
Note that $u(t)$ is decreasing. Using H\"{o}lder's inequality (recall $q \leq 2$) twice, we esimtate:
\begin{eqnarray*}
u(t)^{\frac{q}{q-1}} & \leq & \brac{\int |P_{t/2} f|^q d\mu}^{\frac{q}{q-1}} \leq \brac{\int |P_{t/2} f| d\mu}^{\frac{(2-q)q}{q-1}} \brac{\int |P_{t/2} f|^2 d\mu}^q  \\
& \leq & \norm{f}^{\frac{(2-q)q}{q-1}}_{L_1(\mu)} \brac{\int  f P_t f d\mu}^q \leq 
\norm{f}^{\frac{(2-q)q}{q-1}}_{L_1(\mu)} \norm{f}^q_{N(\mu)^*} \norm{P_t f}^q_{N(\mu)} \\
& \leq &  \norm{f}^{\frac{(2-q)q}{q-1}}_{L_1(\mu)} \norm{f}^q_{N(\mu)^*} \frac{1}{D^q} \int |\nabla P_t f|^q d\mu \\
& \leq &  
\frac{\norm{f}^{\frac{(2-q)q}{q-1}}_{L_1(\mu)} \norm{f}^q_{N(\mu)^*}}{D^q} 
\brac{\int |P_t f|^q d\mu}^{\frac{2-q}{2}} \brac{ \int  |P_t f|^{q-2} | \nabla P_t f |^2 d\mu }^{\frac{q}{2}} ~.
\end{eqnarray*}
Rearranging terms, we see that:
\[
 u'(t) \leq -q(q-1) \frac{D^2}{\norm{f}^{\frac{2(2-q)}{q-1}}_{L_1(\mu)} \norm{f}^2_{N(\mu)^*}} u(t)^{1 + \frac{2}{q(q-1)}} ~.
\]
Setting $v(t) = u(t)^{-\frac{2}{q(q-1)}}$, we obtain:
\[
 v'(t) \geq \frac{2 D^2}{\norm{f}^{\frac{2(2-q)}{q-1}}_{L_1(\mu)} \norm{f}^2_{N(\mu)^*}} ~,
\]
and the assertion follows after integrating in $t$.
\end{proof}

\begin{proof}[Proof of Theorem \ref{thm:Orlicz-implies-I} for $q \leq 2$ - semi-group approach]
We begin as in the proof of the $q \geq 2$ case above, obtaining for a Borel set $A \subset \Omega$ and any time $0 \leq t \leq 1/(2\kappa)$:
\[
2 \sqrt{t} \mu^+(A) \geq \int \abs{\chi_A - P_t \chi_A} d\mu.
\]
Denote $f = \chi_A - \mu(A)$. Since $|x|^q - |y|^q \leq q |x - y|$ for all $|x|,|y| < 1$, we have:
\[
2 \sqrt{t} \mu^+(A) \geq \int \abs{f - P_t f} d\mu \geq \frac{1}{q} \brac{\int |f|^q d\mu - \int |P_t f|^q d\mu} ~,
\]
and using Proposition \ref{prop:new2}, we obtain for $0 \leq t \leq 1/(2\kappa)$:
\begin{equation} \label{eq:optimize-on-2}
 2 \sqrt{t} \mu^+(A) \geq \frac{1}{q} \int |f|^q d\mu \brac{ 1 - \brac{1 + 2 M t }^{-\frac{q(q-1)}{2}}  } ~, 
\end{equation}
where:
\[
M := \frac{D^2 (\int|f|^q d\mu)^{\frac{2}{q(q-1)}}}{\norm{f}^{\frac{2(2-q)}{q-1}}_{L_1(\mu)} \norm{f}^2_{N(\mu)^*}} ~.
\]
As in the proof of the $q \geq 2$ case, it is easy to verify that:
\[
 M \geq \frac{D^2 \brac{\mu(A)(1-\mu(A))}^{\frac{2}{q(q-1)}} N^\wedge(\min(\mu(A),1-\mu(A)))^2 } { \brac{2 \mu(A) (1-\mu(A))}^{\frac{2(2-q)}{q-1} + 2}} ~,
\]
which simplifies to:
\[
 M \geq E:= \frac{D^2 N^\wedge(\min(\mu(A),1-\mu(A)))^2 }{ 2^{\frac{2}{q-1}} \brac{\mu(A)(1-\mu(A))}^{\frac{2}{q}} } ~.
\]

As usual, we will need to optimize (\ref{eq:optimize-on-2}) in $t$. Set $t_0 := 1/E$.
\begin{enumerate}
 \item 
If $t_0 \leq 1/(2\kappa)$, we obtain:
\begin{eqnarray*}
 \mu^+(A) &\geq& \frac{c_1(q-1)}{\sqrt{t_0}} \mu(A)(1-\mu(A)) \\
&\geq& \frac{c_2 (q-1)}{2^{\frac{1}{q-1}}} D \min(\mu(A),1-\mu(A))^{1-1/q} N^\wedge(\min(\mu(A),1-\mu(A))) ~,
\end{eqnarray*}
for some numeric constants $c_1,c_2>0$. 
\item If $t_0 > 1/(2\kappa)$, we evaluate (\ref{eq:optimize-on-2}) at $t_1 = 1/(2\kappa)$. Since $E t_1 < 1$, we obtain:
\begin{eqnarray*}
 \mu^+(A) &\geq& \frac{c_3(q-1) E t_1 }{\sqrt{t_1}} \mu(A)(1-\mu(A)) \geq  \\
&\geq& \frac{c_4 (q-1)}{4^{\frac{1}{q-1}}} \frac{D^2}{\sqrt{\kappa}} \min(\mu(A),1-\mu(A))^{1-2/q} N^\wedge(\min(\mu(A),1-\mu(A)))^2 ~,
\end{eqnarray*}
where $c_3,c_4 > 0$ are numeric constants. Using that $N^\wedge(t)/t^{1/q}$ is non-increasing, we conclude that:
\[
 \mu^+(A) \geq \frac{c_5 (q-1) N^\wedge(1/2)}{4^{\frac{1}{q-1}}} \frac{D^2}{\sqrt{\kappa}} \min(\mu(A),1-\mu(A))^{1-1/q} N^\wedge(\min(\mu(A),1-\mu(A))) ~.
\]
\end{enumerate}
Combining both cases, Theorem \ref{thm:Orlicz-implies-I} follows for $1 < q \leq 2$. 
\end{proof}

\subsubsection{Capacity Approach} \label{subsubsec:capacity-approach}

To complete a circle as we conclude this work, we present a second proof using capacities following \cite[Theorem 4.5]{EMilmanRoleOfConvexityInFunctionalInqs}.

\begin{proof}[Proof of Theorem \ref{thm:Orlicz-implies-I} for $q \leq 2$ - capacity approach]

The assumption (\ref{eq:Orlicz-inq}) implies by Proposition \ref{prop:Capq-Lq-weak} that:
\[
 Cap_q(t,1/2) \geq D N^\wedge(t) \;\;\; \forall t \in [0,1/2] ~,
\]
where $N^\wedge(t) / t^{1/q}$ is non-increasing by our assumption on $N$. 
Using Proposition \ref{prop:increase-cap-q} (with $q_0=q,q=2$)
to pass from $Cap_q$ to $Cap_2$, we obtain that:
\[
 Cap_2(t,1/2) \geq \frac{(1-\frac{2}{q^*})^{1/2}}{(\frac{q^*}{2} - 1)^{1/q^*}} D \brac{\int_t^{1/2} \frac{ds}{(s-t)^{2/q^*} N^\wedge(s)^2}}^{-1/2} \;\;\; \forall t \in [0,1/2] ~.
\]
Next, we modify $N^\wedge(t)$ when $t \geq 1/2$ as follows:
\[
 N_0^\wedge(t) = \begin{cases} N^\wedge(t) & t \in [0,1/2] \\ N^\wedge(1/2) 2^{1/q} t^{1/q} & t \in [1/2,\infty) \end{cases} ~,
\]
so that $N_0^\wedge(t)/t^{1/q}$ is still non-increasing. Using \cite[Lemma 4.2 and Remark 4.3]{EMilmanRoleOfConvexityInFunctionalInqs}, it follows that there exists a numeric constant $c_1>0$ so that:
\[
 Cap_2(t,1/2) \geq c_1 D N_2^\wedge(t) \;\;\; \forall t \in [0,1/2] ~,
\]
where $N_2 \in \J$ is a function so that:
\[
 N_2^\wedge(t) := \frac{1}{\brac{\int_t^{\infty} \frac{ds}{s^{2/q^*} N^\wedge_0(s)^2}}^{1/2}} ~.
\]
Moreover, by \cite[Lemma 4.4]{EMilmanRoleOfConvexityInFunctionalInqs}, $N_2$ is in fact a convex function and $N_2(t)^{1/2}/t$ is non-decreasing.
Proposition \ref{prop:Capq-Lq} then implies that:
\[
 \forall f \in \F \;\;\;\; \frac{c_1}{4} D \norm{f - M_\mu f}_{N_2(\mu)} \leq \norm{ \abs{\nabla f}
}_{L_2(\mu)} ~.
\]
We can now apply the case $q=2$ of Theorem \ref{thm:Orlicz-implies-I}, and conclude that:
\[
 \tilde{I}(t) \geq \min(c_2,N_2^\wedge(1/2)) \min\brac{D,\frac{D^2}{\sqrt{\kappa}}} t^{1/2} N_2^\wedge(t) \;\;\; \forall t \in [0,1/2] ~.
\]
with $c_2>0$ a numeric constant. Defining:
\begin{equation} \label{eq:CNq}
 C_{N,q} := \min(c_2,N_2^\wedge(1/2))  \inf_{0<t<1/2} \frac{t^{1/q - 1/2}}{\brac{\int_t^{\infty}
\frac{N^\wedge_0(t)^2 ds}{s^{2/q^*} N^\wedge_0(s)^2  }}^{1/2}} ~,
\end{equation} 
since $N^\wedge_0(t) = N^\wedge(t)$ for $t\in[0,1/2]$, this implies:
\[
 \tilde{I}(t) \geq C_{N,q} \min\brac{D,\frac{D^2}{\sqrt{\kappa}}} t^{1-1/q} N^\wedge(t) \;\;\; \forall t \in [0,1/2],
\]
as required. This concludes the proof under the additional assumption that $C_{N,q} > 0$. 
\end{proof}

\begin{rem}
Note the similarity between the definitions of $C_{N,q}$ in (\ref{eq:CNq}) and $B_{N,q}$ in (\ref{eq:B-estimate}). 
\end{rem}

This proof has the advantage that one may obtain estimates which do not degrade to $0$ as $q \rightarrow 1$, if the constant $C_{N,q}$ in (\ref{eq:CNq}) may be controlled. Indeed, this is the case for the family of $q$-log-Sobolev inequalities (\ref{eq:q-log-Sob}) for $q\in[1,2]$, discussed in Section \ref{sec:defs}. As already mentioned, it was shown by Bobkov and Zegarlinski that these inequalities may be put in an equivalent form, given by (\ref{eq:N_q-log-Sob}), corresponding to $(\varphi_q,q)$ Orlicz-Sobolev inequalities, where $\varphi_q = t^q \log (1 + t^q)$. It is not hard to verify (see \cite[Corollary 4.8]{EMilmanRoleOfConvexityInFunctionalInqs}) that $C_{\varphi_q,q} \geq c > 0$ uniformly in $q \in [1,2]$, and so we deduce:

\begin{cor}
Under our $\kappa$-semi-convexity assumptions, the $q$-log-Sobolev inequality:
\[
\forall f \in \F \;\;\;\;  D \brac{\int |f|^q \log |f|^q d\mu - \int
|f|^q d\mu \log(\int |f|^q d\mu)}^{1/q} \leq \norm{\abs{\nabla
f}}_{L_q(\mu)} 
\]
with $1 \leq q \leq 2$, implies the following isoperimetric inequality:
\[
 \tilde{I}(t) \geq c \min\brac{D,\frac{D^2}{\sqrt{\kappa}}} t \log^{1/q} 1/t \;\;\; \forall t \in [0,1/2] ~,
\]
where $c>0$ is a numeric constant (independent of $q$). 
\end{cor}

\setlinespacing{1.0}

\bibliographystyle{plain}

\def\cprime{$'$}

\end{document}